\documentclass[leqno,draft]{article}

\def\im{\mathop{\rm Im}}
\def\re{\mathop{\rm Re}}

\newtheorem{theorem}{Theorem}
\newtheorem{lemma}[theorem]{Lemma}
\newtheorem{proposition}[theorem]{Proposition}
\newtheorem{definition}[theorem]{Definition}
\newtheorem{corollary}[theorem]{Corollary}

\newcommand{\begintheorem}{\addtocounter{equation}{1}\begin{theorem}}
\newcommand{\beginlemma}{\addtocounter{equation}{1}\begin{lemma}}
\newcommand{\beginproposition}{\addtocounter{equation}{1}\begin{proposition}}
\newcommand{\begindefinition}{\addtocounter{equation}{1}\begin{definition}}
\newcommand{\begincorollary}{\addtocounter{equation}{1}\begin{corollary}}

\begin{document}

\title{Elements of harmonic analysis, 2}

\author{Stephen William Semmes 	\\
	Rice University		\\
	Houston, Texas}

\date{}

\maketitle

\tableofcontents

\bigskip

        These informal notes are based on a course given at Rice
University in the spring semester of 2004, and much more information
can be found in the references.

\section{Preliminaries}
\label{preliminaries}
\setcounter{equation}{0}

	Let $n$ be a positive integer, and let ${\bf R}^n$, ${\bf
C}^n$ denote the usual real and complex vector spaces consisting of
$n$-tuples of real and complex numbers, using coordinatewise addition
and scalar multiplication.  Of course any real linear transformation
from ${\bf R}^n$ into itself can be extended to a complex-linear
transformation from ${\bf C}^n$ to itself in a simple way.  Similarly,
a real-linear mapping from ${\bf R}^n$ into ${\bf R}$ or ${\bf C}$ can
be extended to a complex-linear mapping into ${\bf C}$.

	If $a$ is a complex number, then $a$ can be written as $a_1 +
i \, a_2$, where $a_1$, $a_2$ are real numbers.  We call $a_1$, $a_2$
the real and imaginary parts of $a$, which may be written as $\re a$,
$\im a$.  The complex conjugate of $a$ is the complex number defined
by
\begin{equation}
	\overline{a} = a_1 - i \, a_2.
\end{equation}
One can easily verify that
\begin{equation}
	\overline{a + b} = \overline{a} + \overline{b},
		\quad \overline{a \, b} = \overline{a} \, \overline{b}
\end{equation}
for all $a, b \in {\bf C}$.

	The absolute value of a real number $x$ is denoted $|x|$ and
defined to be equal to $x$ when $x \ge 0$ and to $-x$ when $x \le 0$.
This is extended to complex numbers by setting
\begin{equation}
	|a| = \big((\re a)^2 + (\im a)^2\big)^{1/2}
\end{equation}
for $a \in {\bf C}$, which is equivalent to saying that $|a|$ is
a nonnegative real number and that
\begin{equation}
	|a|^2 = a \, \overline{a}.
\end{equation}
The triangle inequality for real or complex numbers states that
\begin{equation}
	|a + b| \le |a| + |b|.
\end{equation}
Also,
\begin{equation}
	|a \, b| = |a| \, |b|
\end{equation}
for all $a, b \in {\bf C}$.

	If $x = (x_1, \ldots, x_n), y = (y_1, \ldots, y_n) \in {\bf R}^n$,
then their inner product is defined by
\begin{equation}
	\langle x, y \rangle = \sum_{j=1}^n x_j \, y_j.
\end{equation}
For $z = (z_1, \ldots, z_n), w = (w_1, \ldots, w_n) \in {\bf C}^n$,
their inner product is defined by
\begin{equation}
	\langle z, w \rangle = \sum_{j=1}^n z_j \, \overline{w_j}.
\end{equation}
Let us also put
\begin{equation}
	z \cdot w = \sum_{j=1}^n z_j \, w_j,
\end{equation}
which is the same as $\langle z, w \rangle$ when $z, w \in {\bf R}^n$.

	The standard Euclidean norm on ${\bf R}^n$, ${\bf C}^n$
can be written as
\begin{equation}
	|w| = \biggl(\sum_{j=1}^n |w_j|^2 \biggr)^{1/2}.
\end{equation}
This is equivalent to saying that $|w|$ is a nonnegative real number
and that
\begin{equation}
	|w|^2 = \langle w, w \rangle.
\end{equation}
The Cauchy--Schwarz inequality states that
\begin{equation}
	|\langle z, w \rangle| \le |z| \, |w|
\end{equation}
for all $z, w \in {\bf C}^n$, and in particular when $z, w \in {\bf R}^n$,
and it can be shown using the fact that
\begin{equation}
	\langle z + s \, w, z + s \, w \rangle
\end{equation}
is a nonnegative real number for all scalars $s$.  The triangle
inequality states that
\begin{equation}
	|z + w| \le |z| + |w|
\end{equation}
for all $z, w \in {\bf C}^n$, and it can be shown using the
Cauchy--Schwarz inequality.

	By a multi-index we mean an $n$-tuple $\alpha = (\alpha_1,
\ldots, \alpha_n)$ of nonnegative integers.  The degree of $\alpha$
is defined by
\begin{equation}
	\delta(\alpha) = \sum_{j=1}^n \alpha_j.
\end{equation}
We also get the associated monomials on ${\bf R}^n$, ${\bf C}^n$
given by
\begin{equation}
	x^\alpha = x_1^{\alpha_1} \cdots x_n^{\alpha_n},
		\quad z^\alpha = z_1^{\alpha_1} \cdots z_n^{\alpha_n},
\end{equation}
i.e., the product of the corresponding powers of $x_j$ or $z_j$, as
apprpriate.  If $\alpha_j = 0$ for some $j$, then we interpret that
factor as being equal to $1$.

	A polynomial on ${\bf R}^n$ or ${\bf C}^n$ means a complex
linear combination of monomials.  In particular, a polynomial on ${\bf
R}^n$ extends to a polynomial on ${\bf C}^n$ in a natural way.  To be
precise, by a polynomial on ${\bf C}^n$ we mean a complex polynomial
or holomorphic polynomial, as opposed to a real polynomial which might
also use complex conjugates.

	Suppose that $f(x)$ is a smooth complex-valued function on
${\bf R}^n$.  We write
\begin{equation}
	\frac{\partial}{\partial x_j} f(x)
\end{equation}
for the usual partial derivative of $f$ in the $j$th direction.
If $\alpha$ is a multi-index, then we put
\begin{equation}
	\frac{\partial^{\delta(\alpha)}}{\partial x^\alpha} f(x)
		= \frac{\partial^{\alpha_1}}{\partial x_1^{\alpha_1}}
	  \cdots \frac{\partial^{\alpha_n}}{\partial x_n^{\alpha_n}} f(x).
\end{equation}
In other words, this is the mixed partial derivative of $f$ of order
equal to the degree of $\alpha$ with $\alpha_j$ derivatives in the
$j$th direction for each $j$.

	Now suppose that that $f(z)$ is a smooth complex-valued
function on ${\bf C}^n$, and let us write $x_j$, $y_j$ for the real
and imaginary parts of $z_j$.  We can think of $f(z)$ as a smooth
function of the $2n$ real variables $x_j$, $y_j$, and define partial
derivatives accordingly.  We define complex first-order partial
derivatives by
\begin{equation}
	\frac{\partial}{\partial z_j} f(z)
		= \frac{1}{2} \biggl(\frac{\partial}{\partial x_j}
			- i \frac{\partial}{\partial y_j} \biggr) f(z)
\end{equation}
and
\begin{equation}
	\frac{\partial}{\partial \overline{z_j}}
		= \frac{1}{2} \biggl(\frac{\partial}{\partial x_j}
			+ i \frac{\partial}{\partial y_j} \biggr) f(z).
\end{equation}
We say that $f(z)$ is holomorphic if
\begin{equation}
	\frac{\partial}{\partial \overline{z_j}} f(z) = 0
\end{equation}
for $j = 1, \ldots, n$.

	Constants and complex-linear functions on ${\bf C}^n$
are clearly holomorphic.  The sum and product of two holomorphic
functions on ${\bf C}^n$ are again holomorphic.  As a result,
(complex) polynomials on ${\bf C}^n$ are holomorphic.

\section{Nice functions}
\label{nice functions}
\setcounter{equation}{0}

	Recall that a linear transformation $A : {\bf R}^n \to {\bf
R}^n$ is said to be symmetric if
\begin{equation}
	\langle A(x), y \rangle = \langle x, A(y) \rangle
\end{equation}
for all $x, y \in {\bf R}^n$.  If also
\begin{equation}
	\langle A(x), x \rangle > 0
\end{equation}
when $x \ne 0$, then $A$ is said to be positive-definite.  The
principal axis theorem states that if $A$ is a symmetric linear
transformation on ${\bf R}^n$, then there is an orthonormal basis for
${\bf R}^n$ consisting of eigenvectors for $A$.  Of course the
corresponding eigenvalues are positive real numbers if $A$ is also
positive-definite.

	If $z$ is a complex number, then the exponential of $z$
is defined by
\begin{equation}
	\exp z = \sum_{j=0}^\infty \frac{z^j}{j!}.
\end{equation}
As usual $z^j$ is interpreted as being equal to $1$ when $j = 0$, and
$j!$ is $j$ factorial, the product of the positive integers from $1$
to $j$, which is also interpreted as being equal to $1$ when $j = 0$.
By standard results this power series converges absolutely for all $z
\in {\bf C}$, and converges uniformly on bounded subsets of ${\bf C}$.

	One can check that
\begin{equation}
	\exp (z + w) = \exp(z) \, \exp(w)
\end{equation}
for all complex numbers $z$, $w$.  Specifically, one can formally
multiply the two series on the right and rearrange terms to get the
series on the left, using the binomial theorem, and absolute
convergence of the series is adequate to ensure that the formal
computation is correct.  As a special case,
\begin{equation}
	\exp (z) \exp(-z) = 1
\end{equation}
for all complex numbers $z$, and in particular $\exp(z) \ne 0$
for all $z$.

	One can also check that
\begin{equation}
	\overline{\exp(z)} = \exp (\overline{z}),
\end{equation}
i.e., the complex conjugate of the exponential is the exponential
of the complex conjugate.  If $x$ is a real number, $\exp (x)$
is a real number, and in fact $\exp(x) > 0$ because this is clear
when $x \ge 0$ and when $x < 0$ we have that $\exp (x)$ is the
reciprocal of $\exp(-x) > 0$.  Furthermore,
\begin{equation}
	|\exp (z)| = \exp (\re z)
\end{equation}
for all complex numbers $z$.

	Let us say that a complex-valued function on ${\bf R}^n$ is a
nice function if it is a complex linear combination of functions of
the form
\begin{equation}
	x^\alpha \exp (A(x) \cdot x + b \cdot x),
\end{equation}
where $\alpha$ is a multi-index, $A$ is a positive-definite symmetric
linear transformation on ${\bf R}^n$, and $b \in {\bf C}^n$.
In fact a nice function automatically has a holomorphic extension to
${\bf C}^n$, just using the same formula.

\section{Fourier transforms}
\label{Fourier transforms}
\setcounter{equation}{0}

	If $f(x)$ is a nice function on ${\bf R}^n$, then we can
define the Fourier transform of $f$ by
\begin{equation}
	\widehat{f}(\xi) = 
		\int_{{\bf R}^n} f(x) \, \exp (- 2 \pi i x \cdot \xi) \, dx.
\end{equation}
For $\xi \in {\bf R}^n$ we have that
\begin{equation}
	|\widehat{f}(\xi)| \le \int_{{\bf R}^n} |f(x)| \, dx.
\end{equation}
The Fourier transform $\widehat{f}(\xi)$ of a nice function $f(x)$
makes sense for all $\xi \in {\bf C}^n$, and defines a holomorphic
function on ${\bf C}^n$, because of the rapid decay of nice functions
as $|x| \to \infty$, and one might refer to this holomorphic function
on ${\bf C}^n$ as the Fourier--Laplace transform of $f$.  When $x, \xi
\in {\bf R}^n$ we have that
\begin{equation}
	|\exp (- 2 \pi i x \cdot \xi)| = 1,
\end{equation}
and this is not true in general when $\xi \in {\bf C}^n$, so that we
do not have the same kind of uniform bound for $|\widehat{f}(\xi)|$
for $\xi \in {\bf C}^n$ as for $\xi \in {\bf R}^n$.

	Suppose that $n = 1$, and consider the nice function
\begin{equation}
	f(x) = \exp (- \pi x^2).
\end{equation}
It is well known that
\begin{equation}
	\int_{\bf R} \exp (-\pi x^2) \, dx = 1.
\end{equation}
Clearly the integral is a positive real number, and to determine its
value one can use the well-known trick that the square of the integral
is equal to the analogous integral over ${\bf R}^n$, and that can be
converted to an elementary $1$-dimensional integral via polar
coordinates.

	If $z$ is a real number, then
\begin{eqnarray}
	\int_{\bf R} \exp (-\pi x^2 - 2 \pi x \, z) \, dx
	   & = & \exp (\pi z^2) \int_{\bf R} \exp (-\pi (x + z)^2) \, dx  \\
	   & = & \exp (\pi z^2) \int_{\bf R} \exp (-\pi x^2) \, dx
							\nonumber 	\\
	   & = & \exp (\pi z^2).			\nonumber
\end{eqnarray}
Actually, the conclusion
\begin{equation}
	\int_{\bf R} \exp (-\pi x^2 - 2 \pi x \, z) \, dx
		= \exp (\pi z^2)
\end{equation}
also works if $z$ is a complex number, because both sides of the
equation define holomorphic functions of $z$ which are equal
when $z$ is real.  As a result, if $f(x) = \exp (-\pi x^2)$,
then its Fourier transform is given by
\begin{equation}
	\widehat{f}(\xi) = \exp(- \pi \xi^2),
\end{equation}
and this works for all complex numbers $\xi$.  If we take $f(x)$ to be
the nice function on ${\bf R}^n$ given by $f(x) = \exp (-\pi x \cdot
x)$, then we have that
\begin{equation}
	\widehat{f}(\xi) = \exp (-\pi \xi \cdot \xi)
\end{equation}
for all $\xi \in {\bf C}^n$, because the $n$-dimensional integral in
the Fourier transform reduces to a product of $n$ $1$-dimensional
integrals just computed.

	Suppose that $T$ is a linear transformation on ${\bf R}^n$,
which we can also view as a complex linear transformation on ${\bf
C}^n$ which takes ${\bf R}^n$ to itself.  There is a unique linear
transformation $T^t$ on ${\bf C}^n$, called the transpose of $T$,
which maps ${\bf R}^n$ to itself and satisfies
\begin{equation}
	T(v) \cdot w = v \cdot T^t(w)
\end{equation}
for all $v, w \in {\bf C}^n$.  As usual the standard matrix associated
to $T^t$ is the transpose of the standard matrix associated to $T$,
which is to say that one interchanges the indices of the matrix.

	If $T$ is an invertible linear transformation on ${\bf R}^n$
and $f(x)$ is a nice function on ${\bf R}^n$, then the composition
$f(T(x))$ is also a nice function on ${\bf R}^n$, as one can verify
directly from the definition.  Making a change of variables in the
integral defining the Fourier transform one can also check that the
Fourier transform of $f(T(x))$ is equal to
\begin{equation}
	|\det T|^{-1} \widehat{f}(\widetilde{T}(\xi)),
\end{equation}
where $\det T$ denotes the determinant of $T$ and $\widetilde{T}$ is
the inverse of the transpose of $T$, which is the same as the
transpose of the inverse of $T$.  Note that $T$ is an orthogonal
transformation on ${\bf R}^n$ if and only if $\widetilde{T} = T$, in
which case $|\det T| = 1$ and the Fourier transform of $f(T(x))$ is
equal to $\widehat{f}(T(\xi))$.

	Suppose that $A$ is a positive-definite symmetric linear
transformation on ${\bf R}^n$, and that
\begin{equation}
	f(x) = \exp (- \pi A(x) \cdot x).
\end{equation}
In this case the Fourier transform of $f$ is given by
\begin{equation}
	\widehat{f}(\xi) = (\det A)^{-1/2} \exp (- \pi A^{-1}(\xi) \cdot \xi).
\end{equation}
To see this one can use the fact that $A$ is diagonalized in an
orthonormal basis to reduce the $n$-dimensional integral to a product
of $n$ $1$-dimensional integrals.  Alternatively, one can derive
this from the change of variables formula discussed in the
previous paragraph, with $A = T^t \circ T$.

	Let $f(x)$ be a nice function on ${\bf R}^n$, and let $b$ be
an element of ${\bf R}^n$.  Thus $f(x) \exp (- 2 \pi i x \cdot b)$ is
also a nice function on ${\bf R}^n$, and it is easy to see that the
Fourier transform of this function is equal to $\widehat{f}(\xi + b)$.
This extends to $b \in {\bf C}^n$, because both functions are
holomorphic in $b$ and agree when $b \in {\bf R}^n$.

	Let $f(x)$ be a nice function on ${\bf R}^n$, and let $\alpha$
be a multi-index.  By differentiating under the integral sign we
have that
\begin{equation}
  \frac{\partial^{\delta(\alpha)}}{\partial \xi^\alpha} \widehat{f}(\xi)
\end{equation}
is equal to the Fourier transform of
\begin{equation}
	(-2 \pi i)^{\delta(\alpha)} x^\alpha f(x),
\end{equation}
which is also a nice function on ${\bf R}^n$.  To be a bit more
precise, one can think of $\widehat{f}(\xi)$ as a function on ${\bf R}^n$,
in which case we are taking real derivatives of it, or as a holomorphic
function on ${\bf C}^n$, in which case we are taking complex derivatives.

	By combining these various properties of the Fourier transform
of a nice function we see that the Fourier transforms of the building
blocks for nice functions can be computed explicitly, and that the
Fourier transform of a nice function is a nice function.

\section{Fourier transforms, continued}
\label{Fourier transforms, continued}
\setcounter{equation}{0}

	Nice functions on ${\bf R}^n$ are defined as linear
combinations of the building blocks
\begin{equation}
	x^\alpha \exp (- A(x) \cdot x + b \cdot x),
\end{equation}
where $\alpha$ is a multi-index, $A$ is a positive-definite symmetric
linear transformation on ${\bf R}^n$, and $b \in {\bf C}^n$.
One might instead consider the building blocks
\begin{equation}
	\frac{\partial^{\delta(\alpha)}}{\partial x^\alpha}
		\exp(- A(x) \cdot x + b \cdot x),
\end{equation}
where $\alpha$, $A$, and $b$ are as in the previous case.  It is easy
to see that every building block of the second type can be expressed
as a linear combination of building blocks of the first type, simply
by direct computation.

	Conversely, every building block of the first type can be
expressed as a linear combination of building blocks of the second
type.  One can approach this by approximating a building block of the
first type using derivatives of $\exp (- A(x) \cdot x + b \cdot x)$,
with error terms given by this exponential times monomials of lower
degree.  The two types of building blocks agree when $\alpha = 0$, and
using induction on the degree of $\alpha$ one can show that building
blocks of the first type can be expressed as linear combinations of
building blocks of the second type.

	Thus nice functions on ${\bf R}^n$ can be defined either as
linear combinations of the first kind of building block or of the
second kind of building block, and in particular derivatives of nice
functions are nice functions.  Let $f(x)$ be a nice function on ${\bf
R}^n$, let $\alpha$ be a multi-index, and consider
\begin{equation}
	\frac{\partial^{\delta(\alpha)}}{\partial x^\alpha} f(x).
\end{equation}
The Fourier transform of this nice function is given by
\begin{equation}
	(2 \pi i)^{\delta(\alpha)} \xi^\alpha \, \widehat{f}(\xi),
\end{equation}
by integration by parts.

	If $f(x)$ is a nice function on ${\bf R}^n$ and $a \in {\bf
R}^n$, then $f(x - a)$ is a nice function on ${\bf R}^n$ too.
The Fourier transform of this nice function is equal to
\begin{equation}
	\exp (- 2 \pi i \xi \cdot a) \widehat{f}(\xi),
\end{equation}
by a simple change of variables.  This also works for $a \in {\bf
C}^n$, using the holomorphic extensions of $f$, $\widehat{f}$.

	A corollary of the equivalence of the two kinds of building
blocks for nice functions is that the Fourier transform is a linear
mapping from the vector space of nice functions on ${\bf R}^n$
\emph{onto} the vector space of nice functions on ${\bf R}^n$.
If $\phi(\xi)$ is a nice function on ${\bf R}^n$, consider
\begin{equation}
	\check{\phi}(x) = 
		\int_{{\bf R}^n} \phi(\xi) \, \exp(2 \pi i \xi \cdot x) \, dx,
\end{equation}
which makes sense as a holomorphic function on ${\bf C}^n$, and as a
function on ${\bf R}^n$ in particular.  This transform behaves in
practically the same way as the Fourier transform, and it takes nice
functions to nice functions in particular.  One can check that this is
the inverse of the Fourier transform, which is to say that
$\check{\phi} = f$ when $\phi = \widehat{f}$ and conversely.  In this
regard note that both transforms take $\exp (- \pi z \cdot z)$ to
itself, and that the computation of either transform applied to a nice
function can be reduced to this case using linearity and the various
properties of the transforms.

	Suppose that $f(x)$, $\phi(\xi)$ are nice functions on
${\bf R}^n$.  It is easy to see that
\begin{equation}
	\int_{{\bf R}^n} \phi(\xi) \, \overline{\widehat{f}(\xi)} \, d\xi
   =  \int_{{\bf R}^n} \check{\phi}(x) \, \overline{f(x)} \, dx.
\end{equation}
Specifically, both sides of the equation can be expanded into double
integrals.  The two double integrals are equal by interchanging the
order of integration.

	Because the two transforms are inverses of each other,
it follows that
\begin{equation}
  \int_{{\bf R}^n} \widehat{f}_1(\xi) \, \overline{\widehat{f}_2(\xi)} \, d\xi
	= \int_{{\bf R}^n} f_1(x) \, \overline{f_2(x)} \, dx
\end{equation}
for all nice functions $f_1$, $f_2$ on ${\bf R}^n$.  This is known as
Plancherel's theorem.  Specifically, it follows from the previous
identity with $\phi = \widehat{f}_1$ and $f = f_2$.  As a result
\begin{equation}
	\int_{{\bf R}^n} |\widehat{f}(\xi)|^2 \, d\xi
		= \int_{{\bf R}^n} |f(x)|^2 \, dx
\end{equation}
for all nice functions $f$ on ${\bf R}^n$.

	If $f_1$, $f_2$ are nice functions on ${\bf R}^n$, then the
convolution of $f_1$ and $f_2$ is the function on ${\bf R}^n$ defined
by
\begin{equation}
	(f_1 * f_2)(x) = \int_{{\bf R}^n} f_1(y) \, f_2(x - y) \, dy.
\end{equation}
Clearly this is linear in $f_1$, $f_2$, and one can check that the
convolution of two nice functions is again a nice function.
Moreover,
\begin{equation}
	f_1 * f_2 = f_2 * f_1
\end{equation}
and
\begin{equation}
	(f_1 * f_2) * f_3 = f_1 * (f_2 * f_3)
\end{equation}
when $f_1$, $f_2$, $f_3$ are nice functions on ${\bf R}^n$.

	The convolution of two nice functions $f_1$, $f_2$ on ${\bf
R}^n$ is also characterized by the property that the Fourier transform
of $f_1 * f_2$ is equal to the product of the Fourier transforms of
$f_1$, $f_2$.  Similarly, if $\phi_1$, $\phi_2$ are nice functions on
${\bf R}^n$ and $\phi = \phi_1 * \phi_2$, then $\check{\phi} =
\check{\phi}_1 \, \check{\phi}_2$.

\end{document}